\documentclass[12pt]{amsart}
\usepackage[left=3cm, right=3cm, top=2.5cm, bottom=2.5cm]{geometry}

\usepackage{amscd}
\usepackage{amsmath, amssymb, comment}
\usepackage{amsfonts}
\usepackage{enumerate}
\usepackage{color}
\usepackage{mathtools}

\usepackage{hyperref}
\usepackage{url}
\newcommand{\de}{\partial}

\newcommand{\ov}[1]{\overline{#1}}

\newcommand{\ddbar}{\sqrt{-1} \partial \overline{\partial}}
\newcommand{\Ric}{\mathrm{Ric}}

\newcommand{\tr}{\mathrm{tr}}

\newcommand{\vp}{\varphi}

\newcommand{\ve}{\varepsilon}

\renewcommand{\leq}{\leqslant}
\renewcommand{\geq}{\geqslant}
\renewcommand{\le}{\leqslant}

\newcommand{\be}{\begin{equation}}
\newcommand{\ee}{\end{equation}}

\newcommand{\PSH}{\mathrm{PSH}}

\newcommand{\C}{\mathbb{C}}

\renewcommand{\exp}{\mathrm{exp}}

\newcommand{\CC}{\mathbb {C}}

 \newcommand{\eps}{\varepsilon}

\begin{document}
%\newcounter{remark}
\newcounter{theor}
\setcounter{theor}{1}
\newtheorem{claim}{Claim}
\newtheorem{theorem}{Theorem}[section]
\newtheorem{lemma}[theorem]{Lemma}
\newtheorem{corollary}[theorem]{Corollary}
\newtheorem{proposition}[theorem]{Proposition}
\newtheorem{question}{Question}[section]
\newtheorem{goal}{Goal}[section]
\newtheorem{definition}[theorem]{Definition}
\newtheorem{remark}[theorem]{Remark}
\newtheorem{problem}[theorem]{Problem}

\numberwithin{equation}{section}

\title{Quadratic growth rigidity for the complex Monge--Amp\`ere equation}

\author[J. Chu]{Jianchun Chu}
\address{School of Mathematical Sciences, Peking University, Yiheyuan Road 5, Beijing 100871, People's Republic of China}
\email{jianchunchu@math.pku.edu.cn}

\author[Y. Liu]{Yaxiong Liu}
\address{Mathematical Science Research Center, Chongqing University of Technology, No. 69, Hongguang Avenue, Banan District, Chongqing 400054, People's Republic of China
\newline\hspace*{\parindent}%
Department of Mathematics, University of Maryland,
4176 Campus Dr,
College Park, MD 20742, USA
}
\email{yxliu269@cqut.edu.cn}

\begin{abstract}
For the complex Monge--Amp\`ere equation on $\mathbb{C}^{n}$, we prove the Liouville theorem for solutions with quadratic growth.
\end{abstract}

\subjclass[2020]{Primary 32W20; Secondary 35J96}
\keywords{Complex Monge--Amp\`ere equation, Liouville theorem, Quadratic growth}

\maketitle

\section{Introduction}

On $\C^{n}$, each K\"ahler metric $\omega$ can be written as $\ddbar u$ for some function $u\in\PSH(\C^n)\cap C^{\infty}(\C^{n})$, where $u$ is called the potential function of $\omega$. The associated Ricci form is
\[
\Ric(\omega) = -\ddbar\log\det u_{i\ov{j}}.
\]
This shows deep connection between the Calabi--Yau metric and complex Monge--Amp\`ere equation. In particular, if $u$ satisfies
\begin{equation}\label{CMAE Cn}
\det u_{i\ov{j}} = 1 \ \ \text{on $\C^{n}$},
\end{equation}
then $\omega$ is a Calabi--Yau metric.

In \cite{Tian06}, Tian proved that any complete Calabi--Yau metric on $\mathbb{C}^{2}$ with maximal volume growth is flat, and conjectured that similar result holds on $\mathbb{C}^{n}$ for $n\geq3$. Later, Li \cite{Li19}, Conlon--Rochon \cite{CR21},  Sz\'ekelyhidi \cite{Sze19} constructed the counterexample independently. More precisely, the potential function in these counterexamples satisfies \eqref{CMAE Cn}.
This shows that some additional assumptions are needed for the Liouville theorem. When $\Delta u\leq C$ for some constant $C$, Riebesehl--Schulz \cite{RS84} proved that $u_{i\ov{j}}$ are constant and $\omega=\ddbar u$ is flat. It is very natural to prove the Liouville theorem under some zeroth order assumption. At the AIM workshop ``Nonlinear PDEs in real and complex geometry", Sz\'ekelyhidi \cite[Monge-Amp\`ere equations, Problem 1.4]{STW18} posed the following problem:

\begin{problem}\label{Szekelyhidi problem}
Suppose that $u\in\PSH(\C^n)\cap C^{\infty}(\C^{n})$ is a solution of \eqref{CMAE Cn} and satisfies
\begin{equation}\label{quadratic growth u}
C_{0}^{-1}(|z|^{2}+1) \leq u \leq C_{0}(|z|^{2}+1) \ \ \text{on $\C^{n}$}
\end{equation}
for some constant $C_0>0$, then $u$ is a quadratic polynomial.
\end{problem}

Some process towards this was made recently. In \cite{Wang13}, Wang solved Problem \ref{Szekelyhidi problem} under the assumption $u=|z|^{2}+o(|z|^{2})$ as $|z|\rightarrow\infty$. In \cite{LS21}, Li--Sheng solved Problem \ref{Szekelyhidi problem} under the assumption $\omega=\ddbar u$ is a complete metric on $\mathbb{C}^{n}$.

In this paper, we solved Problem \ref{Szekelyhidi problem} without additional assumptions.
\begin{theorem}\label{main result}
Suppose that $u\in\PSH(\C^n)\cap C^{\infty}(\C^{n})$ is a solution of \eqref{CMAE Cn} and satisfies \eqref{quadratic growth u}, then $u$ is a quadratic polynomial.
\end{theorem}

Let us discuss the proof of Theorem \ref{main result}. The main task is to establish the second order estimate of $u$ in some weak sense. Precisely, we aim to construct a sequence $\{v_{j}\}$ such that $\Delta v_{j}\leq C$ and $v_{j}\to u$ in $L_{\mathrm{loc}}^{1}(\C^{n})$. For the construction of $v_{j}$, motivated by Li--Li--Zhang \cite[(1.3)]{LLZ18}, we consider the solution of Dirichlet problem for complex Monge--Amp\`ere equation with boundary value given by some smooth regularization of $u$ (see \eqref{Dirichlet problem v R t}). In this Dirichlet problem, there are two parameters $R$ and $t$. The key is to choose them suitably such that both the Laplacian estimates and convergence hold (see Proposition \ref{R j t j proposition}). However, the quadratic growth assumption \eqref{quadratic growth u} only gives zeroth order information of $u$. To ``transform" zeroth order information into second order information, the twisted envelope is introduced as an important tool for analysis (see Section \ref{sec:twisted envelopes}). This gives us some quantitative estimate (see Proposition \ref{prop:Delta-P_a(u)-leq} and Lemma \ref{lem:upperbound-E}) which are sufficient, together with some estimates for the Dirichlet problem (see Section \ref{sec:estimates for the Dirichlet problem}), to push the argument through.

\bigskip

The rest of the paper is organized as follows. In Section \ref{sec:estimates for the Dirichlet problem}, we collect some estimates for the Dirichlet problem. In Section \ref{sec:twisted envelopes}, the twisted envelope is introduced and investigated. The proof of Theorem \ref{main result} will be given in Section \ref{sec:proof of main result}.

\bigskip

{\bf Acknowledgments:}
The first-named author was partially supported by National Key R\&D Program of China 2024YFA1014800 and 2023YFA1009900, and NSFC grants 12471052 and 12271008. The second-named author was partially supported by NSFC grant 12601149. Part of this work was carried out while the second-named author was visiting the Institute for Theoretical Sciences at Westlake University, which he would like to thank for the hospitality and support.

\bigskip

{\bf Declaration on the use of AI: } During discussions between the authors and GPT-6 Astra, GPT-6 Astra introduced the notion of the twisted envelope and suggested establishing some quantitative estimates. The proof of Proposition \ref{prop:Delta-P_a(u)-leq} and Lemma \ref{lem:upperbound-E} were critically assisted by GPT-6 Astra.

All mathematical arguments were written by the authors, who take full responsibility for the content of this paper.

\section{Estimates for the Dirichlet problem}\label{sec:estimates for the Dirichlet problem}

Let us focus on the following Dirichlet problem:
\begin{equation}\label{Dirichlet problem v}
\begin{cases}
\ \det v_{i\ov{j}} = 1 & \mbox{in $B_{R}$}, \\[1mm]
\ v = \vp & \mbox{on $\de B_{R}$}.
\end{cases}
\end{equation}
In this section, we collect some second order estimates of $v$ for later use.

\begin{theorem}[Interior estimate, Bedford--Taylor \cite{BT76}]\label{interior estimate}
Suppose that $v\in\PSH(B_{R})\cap C^{\infty}(\ov{B}_{R})$ is a solution of \eqref{Dirichlet problem v}. Then there exists a constant $C$ depending only on $n$ such that
\[
\|\nabla^{2}v\|_{L^{\infty}(B_{R/2})}
\leq C\|\mathrm{Hess}_{\de B_{R}}\vp\|_{L^{\infty}(\de B_{R})}+C.
\]
\end{theorem}

\begin{proposition}[Global estimate]\label{global estimate}
Suppose that $v\in\PSH(B_{R})\cap C^{\infty}(\ov{B}_{R})$ is a solution of \eqref{Dirichlet problem v} and satisfies
\begin{equation}\label{quadratic growth v}
K^{-1}(|z|^{2}+1) \leq v(z) \leq K(|z|^{2}+1) \ \ \text{in $B_{R}$}
\end{equation}
for some constant $K$. Set $M_{R}:=\max_{\de B_{R}}\Delta v+1$. Then there exists a constant $C$ depending only on $K$ and $n$ such that
\[
\Delta v(z) \leq C\exp\left(\frac{C\log M_{R}}{R^{2}+1}\cdot v(z)\right) \ \ \text{in $B_{R}$}.
\]
\end{proposition}

\begin{proof}
For convenience, set
\[
\alpha := \ddbar v, \ \ \ \beta := \ddbar|z|^{2}, \ \ \ A := \frac{R^{2}+1}{K\log M_{R}}.
\]
Then $\Delta v=\tr_{\beta}\alpha$. Consider the quantity
\[
Q := \log\tr_{\beta}\alpha+\frac{|z|^{2}}{v+A}-\frac{v}{A}.
\]
Let $x_{0}$ be the maximum point of $Q$. It suffices to show $Q(x_{0})\leq C$. There are two cases.

\bigskip
\noindent
{\bf Case 1.} $x_{0}\in\de B_{R}$.
\bigskip

Using the definition of $M_{R}$ and \eqref{quadratic growth v}, at $x_{0}$, we have
\[
Q\leq \log M_{R}+K-\frac{R^{2}+1}{KA} \leq K,
\]
as required.

\bigskip
\noindent
{\bf Case 2.} $x_{0}\in B_{R}$.
\bigskip

By the maximum principle, at $x_{0}$, we have
\[
0 \geq \Delta_{\alpha}Q = \Delta_{\alpha}\log\tr_{\beta}\alpha+\Delta_{\alpha}\left(\frac{|z|^{2}}{v+A}\right)-\frac{\Delta_{\alpha}v}{A}.
\]
Since $\alpha$ is Ricci-flat and $\beta$ is flat, then Yau's Schwarz lemma \cite{Yau78} (see also e.g. \cite[Chap 3, Lemma~3.7]{Sze14}) implies $\Delta_{\alpha}\log\tr_{\beta}\alpha\geq 0$. Together with $\Delta_{\alpha}v=n$, we obtain
\begin{equation}\label{global estimate eqn 1}
0 \geq \Delta_{\alpha}\left(\frac{|z|^{2}}{v+A}\right)-\frac{n}{A}.
\end{equation}
Direct calculation shows
\begin{equation}\label{global estimate eqn 2}
\begin{split}
\Delta_{\alpha}\left(\frac{|z|^{2}}{v+A}\right)
= {} & \frac{\Delta_{\alpha}|z|^{2}}{v+A}-\frac{|z|^{2}\Delta_{\alpha}v}{(v+A)^{2}}-\frac{v^{i\ov{j}}\left(v_{i}(|z|^{2})_{\ov{j}}
+(|z|^{2})_{i}v_{\ov{j}}\right)}{(v+A)^{2}}+\frac{2|z|^{2}v^{i\ov{j}}v_{i}v_{\ov{j}}}{(v+A)^{3}} \\
= {} & \frac{\tr_{\alpha}\beta}{v+A}-\frac{n|z|^{2}}{(v+A)^{2}}-\frac{v^{i\ov{j}}(v_{i}z_{j}
+z_{\ov{i}}v_{\ov{j}})}{(v+A)^{2}}+\frac{2|z|^{2}v^{i\ov{j}}v_{i}v_{\ov{j}}}{(v+A)^{3}}.
\end{split}
\end{equation}
We claim that
\begin{equation}\label{global estimate claim}
\Delta_{\alpha}\left(\frac{|z|^{2}}{v+A}\right) \geq \frac{\tr_{\alpha}\beta}{2(v+A)}-\frac{n|z|^{2}}{(v+A)^{2}}.
\end{equation}
Indeed, if $x_{0}=0$, \eqref{global estimate claim} follows from \eqref{global estimate eqn 2} directly. Otherwise, by the Cauchy-Schwarz inequality,
\[
\frac{v^{i\ov{j}}(v_{i}z_{j}+z_{\ov{i}}v_{\ov{j}})}{(v+A)^{2}}
\leq \frac{2|z|^{2}v^{i\ov{j}}v_{i}v_{\ov{j}}}{(v+A)^{3}}+\frac{v^{i\ov{j}}z_{i}z_{\ov{j}}}{2(v+A)|z|^{2}}.
\]
Combining this with \eqref{global estimate eqn 2}, we obtain \eqref{global estimate claim}. Substituting \eqref{global estimate claim} into \eqref{global estimate eqn 1},
\[
0 \geq \frac{\tr_{\alpha}\beta}{2(v+A)}-\frac{n|z|^{2}}{(v+A)^{2}}-\frac{n}{A}.
\]
Together with \eqref{quadratic growth v},
\[
\tr_{\alpha}\beta \leq 2(v+A)\left(\frac{n}{A}+\frac{n|z|^{2}}{(v+A)^{2}}\right)
\leq 2(v+A)\left(\frac{n}{A}+\frac{nKv}{(v+A)^{2}}\right) \leq 2n\left(\frac{v}{A}+1+K\right).
\]
Since $\det(v_{i\ov{j}})=1$, then $\alpha^{n}=\beta^{n}$ and so
\[
\tr_{\beta}\alpha \leq \frac{(\tr_{\alpha}\beta)^{n-1}}{(n-1)!} \cdot \frac{\alpha^{n}}{\beta^{n}}
\leq \frac{(2n)^{n-1}}{(n-1)!}\left(\frac{v}{A}+1+K\right)^{n-1}.
\]
Combining this with \eqref{quadratic growth v} and elementary inequality $(n-1)\log(x+1+K)-x\leq C$ for $x>0$,
\[
Q \leq \log \frac{(\tr_{\alpha}\beta)^{n-1}}{(n-1)!}+(n-1)\log\left(\frac{v}{A}+1+K\right)+\frac{Kv}{v+A}-\frac{v}{A} \leq C, \\
\]
as required.
\end{proof}

\section{Twisted envelopes on $\mathbb{C}^n$}
\label{sec:twisted envelopes}

Let $f$ be a function on $\C^{n}$. The envelope $P(f)$ is defined as
\begin{align*}
	P(f)\coloneqq \big(\sup\{v\leq f: v\in\PSH(\C^{n}) \}\big)^*,
\end{align*}
where $*$ denotes upper semi-continuous regularization. As a generalization, for $a>0$, the twisted envelope is defined as
\begin{align*}
	P_{a}(f)\coloneqq\big(\sup\{v\leq f: v-a|z|^2\in\PSH(\C^{n}) \}\big)^*.
\end{align*}
The following is the key proposition in this section.
\begin{proposition}\label{prop:Delta-P_a(u)-leq}
Suppose that $u\in\PSH(\C^n)\cap C^{\infty}(\C^{n})$ is a solution of \eqref{CMAE Cn} and satisfies \eqref{quadratic growth u}. For any $0<a<(2C_{0})^{-1}$ where $C_{0}$ is the constant in \eqref{quadratic growth u}, we have $0\le\Delta P_a(u)\leq a^{1-n}$ in the distributional sense.
\end{proposition}

We will first prove some envelope convergence lemma and then give the proof of Proposition \ref{prop:Delta-P_a(u)-leq}.

\subsection{Envelope convergence lemma}

\begin{lemma}
\label{lem:envelope-incre-convergence}
	If $f$ is continuous and satisfies
\[
c_{0}|z|^2 + b_{0} \le f(z)\le B_{0}(|z|^2+1) \ \ \text{on $\C^{n}$}
\]
for some $b_0,c_{0},B_{0}>0$, then for every $0<\gamma<c_{0}$ and $M>0$, set
\begin{equation}\label{f M F M F}
f_M\coloneqq\min\{f,\gamma|z|^2+M\}, \ \
F_M\coloneqq P(f_M), \ \
F\coloneqq P(f).
\end{equation}
Then $F_{M}\nearrow F$ in $C_{\mathrm{loc}}^{0}(\C^{n})$ as $M\to\infty$.
\end{lemma}

\begin{proof}
The argument will be divided into three steps.

\bigskip
\noindent
{\bf Step 1.} $F_M$ is continuous.
\bigskip

Define $\vp_M\coloneqq f_M-\gamma|z|^2=\min\{f-\gamma|z|^2,M\}$. Note that $f-\gamma|z|^2\geq b_0+(c_{0}-\gamma)|z|^2$. Then we have $\vp_M(z)=M$ once $|z|^2\geq\frac{M-b_0}{c_{0}-\gamma}$. This shows that, for any $\eps>0$, there exists $\delta>0$ such that for any $z\in\CC^n$ and $|w|<\delta$, $\vp_M(z-w)\leq\vp_M(z)+\eps$. For $|w|<\delta$, define the plurisubharmonic function by
\begin{align*}
		\psi_w(z)\coloneqq F_M(z-w)+\gamma(|z|^2-|z-w|^2)-\eps.
\end{align*}
Since $F_M\leq f_M$, then
\begin{align*}
		\psi_w(z)\leq\gamma|z|^2+\vp_M(z-w)-\eps
		\leq\gamma|z|^2+\vp_M(z)=f_M(z),
\end{align*}
which implies $\psi_w(z)\leq F_M(z)$. For any $z_{0}\in\C^{n}$, we have
\begin{align*}
		F_M(z_0+w) \geq \psi_w(z_0+w) = F_M(z_0)+\gamma(|z_0+w|^2-|z_0|^2)-\eps.
\end{align*}
It then follows that
\begin{align*}
		\liminf_{w\rightarrow0}F_M(z_0+w)\geq F_M(z_0)-\eps = \limsup_{w\rightarrow0}F_M(z_0+w)-\eps,
\end{align*}
where we used the upper-continuity of $F_{M}$ in the equality. Letting $\eps\rightarrow0$, we obtain that $F_M$ is continuous at $z_{0}$.

\bigskip
\noindent
{\bf Step 2.} If $v\in\PSH(\CC^n)$ satisfies $v\leq f$ and $v(z)\leq \sigma|z|^2+C$ for some $\sigma\in(0,c_{0})$ and $C>0$, then $v\leq \tilde F\coloneqq\lim_{M\rightarrow\infty}F_M$.
\bigskip

For any $z_0\in\CC^n$, choose $R$ sufficiently large such that $R>|z_0|$ and $(c_{0}-\sigma)R^2>C-b_0$. Define the plurisubharmonic function $g_{R}$ by
\[
g_R(z)\coloneqq b_0+c_{0}(2R|z|-R^2).
\]
Then $b_0+c_{0}|z|^2-g_R(z)=c_{0}(|z|-R)^2\geq0$ and so
\[
g_R(z) \leq b_0+c_{0}|z|^2 \leq f.
\]
If $z\in\de B_{R}$, then $g_R(z)=b_0+c_{0}R^2>\sigma R^2+C\geq v(z)$. This shows that $\sup_{\de B_{R}}(v-g_{R})<0$. Since $(v-g_{R})$ is upper-continuous, then $v-g_{R}<0$ near $\de B_{R}$. Thus we can construct the plurisubharmonic function
\begin{align*}
\tilde g_R(z)
:=
\begin{cases}
\ \max\{v(z),g_R(z)\}, &|z|<R, \\[1mm]
\ g_R(z), &|z|\geq R.
\end{cases}
\end{align*}
It follows from $g_{R}\leq f$ and $v\leq f$ that $\tilde g_R\leq f$. Since $\tilde g_R$ has only linear growth, then there exists a constant $M\geq b_0$ such that $\tilde g_R\leq\gamma|z|^2+M$ on $\CC^n$. Thus 
\[
v(z_0)\leq\tilde g_R(z_0)\leq F_M(z_0)\leq \tilde F(z_0).
\]

\bigskip
\noindent
{\bf Step 3.} $F_{M}\nearrow F$ in $C_{\mathrm{loc}}^{0}(\C^{n})$ as $M\to\infty$.
\bigskip

Choose $\theta\in(0,1)$ sufficiently small such that $\sigma\coloneqq\theta B_{0}+(1-\theta)\gamma<c_{0}$. For any $M\geq b_0$, we consider the plurisubharmonic function by
\[
G_M\coloneqq\theta F+(1-\theta)F_M,
\]
It is clear that $G_{M}\leq f$ and
\begin{align*}
G_M(z)\leq\theta B_{0}(1+|z|^2)+(1-\theta)(\gamma|z|^2+M)
= \sigma|z|^2+\theta B_{0}+(1-\theta)M.
\end{align*}
Then by Step 2, we have $G_{M}\leq\tilde F$, which implies $\theta F+(1-\theta)F_M\leq\tilde F$. Letting $M\rightarrow\infty$, we obtain $\theta F+(1-\theta)\tilde F\leq\tilde F$ and so $F\leq \tilde F$. Together with $\tilde F\leq F$, we obtain $\lim_{M\rightarrow\infty}F_M=F$. Since $F_{M}$ is increasing with respect to $M$, then by Step 1, $F$ is lower semi-continuous, and then continuous. By Dini's theorem, we obtain the required convergence.	
\end{proof}

\subsection{Proof of Proposition \ref{prop:Delta-P_a(u)-leq}}

By Lemma~\ref{lem:envelope-incre-convergence}, Proposition~\ref{prop:Delta-P_a(u)-leq} follows from Lemma~\ref{lem:Delta-F_M} (see below) immediately. Now we set $f\coloneqq u-a|z|^2$, and still use notations in \eqref{f M F M F}. Recall that $f_M:=\min\{f,\gamma|z|^2+M\}$ and $F_M:=P(f_M)$.
\begin{lemma}
\label{lem:Delta-F_M}
For $0<a<(2C_0)^{-1}$, we have $\Delta F_M\leq a^{1-n}-a$ in the distributional sense.
\end{lemma}

\begin{proof}
Set $b_0:=C_0^{-1}$ and $c_{0}:=C_0^{-1}-a>\frac{1}{2}C_0^{-1}$. Then by \eqref{quadratic growth u}, we see that
\begin{align*}
c_{0}|z|^2+b_0 \le f(z)\leq C_0(|z|^2+1).
\end{align*}

\bigskip
\noindent
{\bf Step 1.} If $\ddbar f(x)\geq0$, then $\Delta f(x)\leq a^{1-n}-a$.
\bigskip

Let $\lambda_1,\ldots,\lambda_{n}$ be the eigenvalues of $\ddbar f(x)$. Then by \eqref{CMAE Cn}, we obtain
\begin{align*}
1 = \det u_{i\ov j}(x) = \prod_{i=1}^n(a+\lambda_i) \geq a^n+a^{n-1}\sum_{i=1}^n\lambda_i,
\end{align*}
as required. As an application, we have
\begin{equation}\label{n gamma a n}
n\gamma\leq a^{1-n}-a.
\end{equation}
Indeed, note that
\begin{align*}
    f-\gamma|z|^2\geq b_0+(c_{0}-\gamma)|z|^2\rightarrow\infty, \ \ \text{ as }|z|\rightarrow\infty.	
\end{align*}
Then it admits a global minimum at some point $z'$, at which $\Delta f(z')\geq n\gamma$. Thus we obtain \eqref{n gamma a n}.

\bigskip
\noindent
{\bf Step 2.} For any $\delta>0$, set $\Lambda\coloneqq \frac{a^{1-n}-a}{n}+\delta$. Then 
\begin{align}
\label{eq:F_M,eps-leq-F_M}
F_M*\rho_\eps-F_M\leq \Lambda m_2\eps^2,
\end{align}
where $\rho_{\ve}$ and $m_{2}$ are defined in \eqref{rho ve}.
\bigskip

Let $\rho\in C^{\infty}(\C^{n})$ be a radial mollifier, that is, $\rho$ is a radial function satisfying
\[
0 \leq \rho \leq 1, \ \ \mathrm{supp}\rho \subset B_{1}, \ \ \int_{B_{1}}\rho = 1.
\]
For $\ve\in(0,1)$, set
\begin{equation}\label{rho ve}
\rho_{\ve}(x) \coloneqq \ve^{-2n}\rho(x/\ve), \ \ m_2\coloneqq\int |y|^2\rho(y)\,\mathrm{d}y.
\end{equation}
Define the smooth plurisubharmonic function by
\begin{align*}
F_{M,\eps}\coloneqq F_M*\rho_\eps-\Lambda m_2\eps^2.
\end{align*}
Since $F_M\leq\gamma|z|^2+M$ and $\gamma\leq \frac{a^{1-n}-a}{n}<\Lambda$ (see \eqref{n gamma a n}), then
\begin{align*}
     	F_M*\rho_\eps\leq (\gamma|z|^2+M)*\rho_\eps
     	=(\gamma|z|^2+M)+\gamma m_2\eps^2
     	\leq (\gamma|z|^2+M)+\Lambda m_2\eps^2,
\end{align*}
which implies $F_{M,\eps}\leq \gamma|z|^2+M$. We claim that $F_{M,\eps}\leq f$. Given this claim, $F_{M,\eps}$ is a admissible function and then $F_{M,\eps}\leq F_M$, which implies
\eqref{eq:F_M,eps-leq-F_M}.

\medskip

We prove the claim by contradiction. Suppose $F_{M,\eps}(z_0)>f(z_0)$ for some $z_0\in\CC^n$. Note that
\begin{align*}
     	F_{M,\eps}-f\leq (\gamma|z|^2+M)-f
     	\leq M-b_0-(c_{0}-\gamma)|z|^2.
\end{align*}
Then $F_{M,\eps}-f$ tends to $-\infty$ as $|z|\rightarrow\infty$ and admits a positive global maximum at some point $z_\eps$. One has
\begin{align*}
     0<F_{M,\eps}(z_\eps)-f(z_\eps)\leq M-b_0-(c_{0}-\gamma)|z_\eps|^2,
\end{align*}
which implies
\begin{align*}
     	|z_\eps|^2\leq\frac{M-b_0}{c_{0}-\gamma}.
\end{align*}
Then $z_\eps\in B_{R_M}$, where $R_M=1+\sqrt{\frac{M-b_0}{c_{0}-\gamma}}$.

By the maximum principle, $\ddbar(F_{M,\eps}-f)(z_\eps)\leq 0$. Since $\ddbar F_{M,\eps}\geq0$, then $\ddbar f(z_\eps)\geq0$.
By Step 1, we have $\Delta f(z_\eps)\leq a^{1-n}-a$.
For any $x\in \overline{B}_{R_M}$, the Taylor expansion implies that
\begin{align*}
        (f*\rho_\eps)(x)-f(x)
      =\frac{m_2}{n}\Delta f(x)\,\eps^2+\eps^2
          r_M(x,\eps),
\end{align*}
where $\sup_{x\in\overline{B}_{R_M}}|r_M(x,\eps)|\rightarrow0$ as $\eps\rightarrow0$. Note that $F_M*\rho_\eps\leq f*\rho_\eps$. Then at $z_\eps$, we have
\[
\begin{split}
0 < {} & F_{M,\eps}(z_\eps)-f(z_\eps) \\[1.5mm]
\leq {} & f*\rho_\eps(z_\eps)-f(z_\eps)-\Lambda m_2\eps^2 \\
\leq {} & \frac{m_2}{n}(a^{1-n}-a)\eps^2+\eps^2r_M(z_\eps,\eps)-\Big(\frac{a^{1-n}-a}{n}+\delta\Big)m_2\eps^2 \\[1mm]
= {} & \eps^2r_M(z_\eps,\eps)-\delta m_2\eps^2.
\end{split}
\]
Choose $\eps>0$ sufficiently small such that $\sup_{x\in\overline{B}_{R_M}}|r_M(x,\eps)|\leq \frac{\delta}{2}m_2$. Then
\begin{align*}
     	0<F_{M,\eps}(z_\eps)-f(z_\eps)
     	\leq -\frac{\delta m_2}{2}\eps^2<0,
\end{align*}
which is a contradiction. This proves the claim.

\bigskip
\noindent
{\bf Step 3.} $0\leq\Delta F_M\leq a^{1-n}-a$ in the distributional sense.
\bigskip

For any $\chi\in C^\infty_c(\C^{n})$ with $\chi\geq0$, integrating \eqref{eq:F_M,eps-leq-F_M} with respect to $\chi\,\mathrm{d}y$ implies
\begin{align*}
     	\int F_M\,\frac{\chi*\rho_\eps-\chi}{\eps^2}\,\mathrm{d}y
     	\leq\Lambda m_2\int\chi\,\mathrm{d}y.
\end{align*}
Again by the Taylor expansion, then $\frac{\chi*\rho_\eps-\chi}{\eps^2}\rightarrow\frac{m_2}{n}\Delta\chi$ uniformly.
Then we obtain
\begin{align*}
     	\frac{m_2}{n}\int F_M\,\Delta\chi \,\mathrm{d}y
     	\leq \left(\frac{a^{1-n}-a}{n}+\delta\right) m_2\int\chi\,\mathrm{d}y.
\end{align*}
Letting $\delta\rightarrow0$, we obtain $\Delta F_M\leq a^{1-n}-a$ in the distributional sense.

\end{proof}

Now, we are in a position to prove Proposition \ref{prop:Delta-P_a(u)-leq}.
\begin{proof}[Proof of \textup{Proposition \ref{prop:Delta-P_a(u)-leq}}]
By Lemma \ref{lem:envelope-incre-convergence}, then $F_M$ locally uniform converges to $P(f)=P_a(u)-a|z|^2$. By Lemma \ref{lem:Delta-F_M}, $0\leq \Delta F_M\leq a^{1-n}-a$. Taking $M\rightarrow\infty$, we obtain $0\leq\Delta P(f)\leq a^{1-a}-a$. Then $\Delta P_a(u)=\Delta P(f)+a\leq a^{1-n}$.
\end{proof}

\section{Proof of Theorem \ref{main result}}\label{sec:proof of main result}
In this section, we will give the proof of Theorem \ref{main result}.

\subsection{Smooth regularization of $u$}
 For $\ve\in(0,1)$, let $\rho_\eps\in C^{\infty}(\C^{n})$ be the scaled mollifier defined by \eqref{rho ve} and set
\[
 u_{\ve} \coloneqq u*\rho_{\ve}.
\]

\begin{lemma}\label{properties of u ve}
The following holds:
\begin{enumerate}\setlength{\itemsep}{1mm}
\item[(i)] $u_{\ve}\geq u$ on $\C^{n}$;
\item[(ii)] $\det(u_{\ve})_{i\ov{j}}\geq1$ on $\C^{n}$;
\item[(iii)] For $R>1$, there is a constant $C$ independent of $R$ such that
\[
\|\mathrm{Hess}_{\de B_{R}}u_{\ve}\|_{L^{\infty}(\de B_{R})} \leq CR^{2}\ve^{-2}+CR\ve^{-1}.
\]
\end{enumerate}
\end{lemma}

\begin{proof}
Since $u$ is subharmonic, then (i) is well-known. For (ii), note that
\[
u_{\ve}(x) = \int_{B_{1}}u(x+\ve y)\rho(y)\,\mathrm{d}y.
\]
Then by the concavity of $\det^{1/n}$ and Jensen's inequality,
\[
\begin{split}
\mathrm{det}^{\frac{1}{n}}(u_{\ve})_{i\ov{j}}(x)
= \mathrm{det}^{\frac{1}{n}}\left(\int_{B_{1}}u_{i\ov{j}}(x+\ve y)\rho(y)\,\mathrm{d}y\right)
\geq \int_{B_{1}}\mathrm{det}^{\frac{1}{n}}(u_{i\ov{j}}(x+\ve y))\rho(y)\,\mathrm{d}y
=  1.
\end{split}
\]
For (iii), it is clear that
\[
u_{\ve}(x) = \ve^{-2n}\int_{B_{\ve}(x)}u(y)\rho\left(\frac{x-y}{\ve}\right)\,\mathrm{d}y.
\]
For $x\in\de B_{R}$ and $y\in B_{\ve}(x)$, using \eqref{quadratic growth u}, we have
\[
|u(y)| \leq C\big((R+\ve)^{2}+1\big) \leq CR^{2}.
\]
It then follows that
\[
|\nabla u_{\ve}(x)| \leq CR^{2}\ve^{-1}, \ \ \
|\nabla^{2} u_{\ve}(x)| \leq CR^{2}\ve^{-2}.
\]
Recall that
\[
\mathrm{Hess}_{\de B_{R}}u_{\ve} = \nabla^{2}u_{\ve}-\langle\nabla u_{\ve},\nu\rangle\cdot\mathrm{II}_{\de B_{R}},
\]
where $\nu$ denotes the unit outer normal of $\de B_{R}$ and $\mathrm{II}_{\de B_{R}}$ denotes the second fundamental form of $\de B_{R}$. Since $|\mathrm{II}_{\de B_{R}}|=1/R$, then
\[
|\mathrm{Hess}_{\de B_{R}}u_{\ve}|(x) \leq CR^{2}\ve^{-2}+CR\ve^{-1}.
\]
\end{proof}

\subsection{Dirichlet problem on $B_{R}$ and smooth approximation of $u$}

Motivated by the construction of Li--Li--Zhang \cite[(1.3)]{LLZ18}, for $R>1$ and $t\in(0,1)$, by \cite[Theorem 1.1]{CKNS85}, the following Dirichlet problem admits a smooth plurisubharmonic solution:
\begin{equation}\label{Dirichlet problem v R t}
\begin{cases}
\ \det(v_{R,t})_{i\ov{j}} = 1 & \mbox{in $B_{R}$}, \\[1mm]
\ v_{R,t} = u_{Rt} & \mbox{on $\de B_{R}$}.
\end{cases}
\end{equation}
Then combining (i), (ii) of Lemma \ref{properties of u ve} with the comparison principle, we obtain
\begin{equation}\label{u u R t v R t}
u \leq u_{Rt} \leq v_{R,t} \ \ \text{on $B_{R}$}.
\end{equation}
Denote the difference of $u$ and $v_{R,t}$ by
\[
D(R,t) := \int_{B_{R}}(v_{R,t}-u).
\]

\begin{lemma}\label{lem:upperbound-D(R,tau)}
For $R>1$ and $t\in(0,1)$, we have
\begin{align*}
		nD(R,t)
		\leq \frac{R}{2}\int_{\partial B_R}(u_{Rt}-u)\cdot\tr(u^{i\ov j}).
	\end{align*}
\end{lemma}
\begin{proof}
	We denote $w_R\coloneqq R^2-|z|^2$ and consider the linearized operator $L_u=u^{i\ov j}\partial_i\ov\partial_{j}$.
	By the AM-GM inequality, one has
	\begin{align*}
		L_u(v_{R,t}-u)=L_u(v_{R,t})-n\geq0.
	\end{align*}
	Since $\det u_{i\ov j}=1$, then
	\begin{align*}
		\sum_i\partial_iu^{i\ov j}=0,  \quad
		\sum_j\partial_{\ov j}u^{i\ov j}=0.
	\end{align*}
	Integration by parts twice yields
	\begin{align*}
		0\leq\int_{B_R}w_RL_u(v_{R,t}-u)
		=&-\int_{B_R}\sum_{i,j}u^{i\ov j}(w_R)_i(v_{R,t}-u)_{\ov j} \\
		=&\int_{B_R}(v_{R,t}-u)L_uw_R
		-\frac{1}{2}\int_{\partial B_R}(v_{R,t}-u)\sum_{i,j}u^{i\ov j}(w_R)_i\nu_j \\
		=&-\int_{B_R}(v_{R,t}-u)\cdot\tr(u^{i\ov j})
		+\frac{1}{2R}\int_{\partial B_R}(v_{R,t}-u)u^{i\ov j}\ov z_i z_j\\
		\leq&-\int_{B_R}(v_{R,t}-u)\cdot\tr(u^{i\ov j})
		+\frac{R}{2}\int_{\partial B_R}(v_{R,t}-u)\cdot\tr(u^{i\ov j}),
	\end{align*}
	where we used $u^{i\ov j}\ov z_i z_j\leq\tr(u^{i\ov j})R^{2}$ on $\de B_{R}$. Combining with $\tr(u^{i\ov j})\geq n$, we obtain the required inequality.
\end{proof}
By Lemma~\ref{lem:upperbound-D(R,tau)} and coarea formula, for any $S,\eps>0$,
\begin{equation}\label{n int S 2S D}
n\int_{S}^{2S}D(R,\eps/R)\,\mathrm{d}R \leq S\int_{\{S<|z|<2S\}}(u_{\eps}-u)\cdot\tr(u^{i\ov j}).
\end{equation}
For simplicity, we denote
$E(S,\eps)\coloneqq \int_{\{S<|z|<2S\}}(u_\eps-u)\cdot\tr(u^{i\ov j})$.

\begin{lemma}
\label{lem:upperbound-E}
	There exists a constant $K>0$ independent of $S$ such that for $\eps=Se^{-nS}$ and $S\gg1$,
	\begin{align*}
		E(S,\eps)\leq KS^{2n+2}(e^{-S}+e^{-(n+1)S}).
	\end{align*}
\end{lemma}
\begin{proof}
For $0<a<(2C_{0})^{-1}$, recall that $P_{a}(u)$ denotes the twisted envelope defined in Section \ref{sec:twisted envelopes}. Since $u$ is continuous, then $P_a(u)\leq u$. The admissible function $(1-aC_0)u+a|z|^2$ gives
\begin{align}
\label{eq:u-P_a(u)-leq}
	0\le u-P_a(u)\le a(C_0u-|z|^2)\leq aC_0u.
\end{align}
Set
\begin{align*}
E_a(S,\eps)
\coloneqq \int_{\{S<|z|<2S\}}\big(P_a(u)*\rho_\eps-P_a(u)\big)\cdot\tr(u^{i\ov j}).
\end{align*}
On the annulus $\{S<|z|<2S\}$, by \eqref{quadratic growth u} and \eqref{eq:u-P_a(u)-leq}, then $0\leq u-P_a(u)\leq Ka S^2$, where $K=K(n,C_0)>0$ is a constant. Note that
\begin{equation}\label{E-E a}
\begin{split}
|E(S,\eps)-E_a(S,\eps)| = {} & \int_{\{S<|z|<2S\}}\big(|u*\rho_\eps-P_a(u)*\rho_\eps|+|u-P_a(u)|\big)\cdot\tr(u^{i\ov j}) \\
= {} & \int_{\{S<|z|<2S\}}\big(|(u-P_a(u))*\rho_\eps|+|u-P_a(u)|\big)\cdot\tr(u^{i\ov j}) \\
\leq {} & KaS^2 \int_{\{S<|z|<2S\}}\tr(u^{i\ov j})
\leq  KaS^{2n+2}.
\end{split}
\end{equation}
Since $\|u\|_{L^{\infty}(B_{2S})}\leq CS^{2}$, then applying integration by parts repeatedly shows
\[
\int_{B_{S}}(\ddbar u)^{n-1}\wedge\ddbar |z|^{2} \leq CS^{2n}.
\]
Together with $\det u_{i\ov{j}}=1$, we obtain
\[
\int_{\{S<|z|<2S\}}\tr(u^{i\ov j}) \leq CS^{2n}.
\]
Combining this with \eqref{E-E a},
\begin{equation}\label{E(S,eps)}
E(S,\eps) \leq E_a(S,\eps)+KaS^{2n+2}.
\end{equation}
On the other hand, Proposition~\ref{prop:Delta-P_a(u)-leq} gives $0\le\Delta P_a(u)\leq a^{1-n}$ in the distributional sense. By the standard subharmonic function theory, we have
	\begin{align*}
		0\leq P_a(u)*\rho_\eps-P_a(u)\leq\frac{a^{1-n}\eps^2}{4n},
	\end{align*}
which implies
\begin{align*}
	E_a(S,\eps)\leq K S^{2n}a^{1-n}\eps^2.
\end{align*}
Combining this with \eqref{E(S,eps)} and choosing $a=e^{-S}$, then we obtain
\begin{align*}
	E(S,\eps)\leq Ke^{-S}S^{2n+2}+KS^{2n}e^{(n-1)S}S^2e^{-2nS}
	\leq KS^{2n+2}(e^{-S}+e^{-(n+1)S}),
\end{align*}
as required.
\end{proof}

\begin{proposition}\label{R j t j proposition}
There exist sequences $\{R_{j}\}$ and $\{t_{j}\}$ such that
\[
R_{j} \to \infty, \ \ t_{j} \to 0, \ \ D(R_{j},t_{j}) \to 0, \ \ \log(t_{j}^{-2}) = o(R_{j}^{2}).
\]
\end{proposition}

\begin{proof}
By Lemma~\ref{lem:upperbound-E}, we can choose a sequence $S_j\rightarrow\infty$ and $\eps_j=S_je^{-nS_j}\rightarrow0$ such that $E(S_j,\eps_j)\rightarrow0$.
By \eqref{n int S 2S D}, there exist $R_j\in[S_j,2S_j]$ and $t_j=\eps_j/R_j$ such that $\log(t_j^{-2})\leq 4nR_j=o(R_j^2)$ and
\begin{align*}
	D(R_j,t_j)\leq\frac{1}{n}E(S_j,\eps_j)\rightarrow0.
\end{align*}
\end{proof}

\subsection{Completion of the proof}

Let $\{R_{j}\}$ and $\{t_{j}\}$ be the sequences obtained in Proposition \ref{R j t j proposition} and set $v_{j}:=v_{R_{j},t_{j}}$. We first show that $v_{j}$ has uniform quadratic growth.

\begin{lemma}\label{uniform quadratic growth}
There is a constant $C$ independent of $j$ such that
\[
C^{-1}(|z|^{2}+1) \leq v_{j}(z) \leq C(|z|^{2}+1) \ \ \text{on $B_{R_{j}}$}.
\]
\end{lemma}

\begin{proof}
The lower bound follows from \eqref{quadratic growth u} and \eqref{u u R t v R t}. For the upper bound, if $|z|<R_{j}^{2}/4$, then the submean inequality shows
\[
v_{j}(z) \leq \frac{1}{|B_{1}(z)|}\int_{B_{1}(z)}(v_{j}-u+u)
\leq \frac{D(R_{j},t_{j})}{|B_{1}(z)|}+\frac{1}{|B_{1}(z)|}\int_{B_{1}(z)}u.
\]
Together with $D(R_{j},t_{j})\to0$ and \eqref{quadratic growth u}, we obtain
\[
v_{j}(z) \leq \frac{C}{|B_{1}|}+\frac{C\big((|z|+1)^{2}+1\big)}{|B_{1}|} \leq C(|z|^{2}+1).
\]
If $|z|\geq R_{j}^{2}/4$, then by the maximum principle, $R_{j}t_{j}\to0$ and \eqref{quadratic growth u},
\[
v_{j}(z) \leq \sup_{\de B_{R_{j}}}u_{R_{j}t_{j}} \leq C\big((R_{j}+C)^{2}+1\big) \leq C(|z|^{2}+1).
\]
\end{proof}

Now we are in a position to prove Theorem \ref{main result}
\begin{proof}[Proof of Theorem \ref{main result}]
Since $D(R_{j},t_{j})\to0$, then $v_{j}\to u$ in $L_{\mathrm{loc}}^{1}(\C^{n})$. To prove Theorem \ref{main result}, it suffices to show that, up to a subsequence, $v_{j}$ converges to a quadratic polynomial pointwisely.

By (iii) of Lemma \ref{properties of u ve} and Theorem \ref{interior estimate}, we have
\begin{equation}\label{M R j 2 estimate}
M_{R_{j}/2}:=\max_{\de B_{R_{j}/2}}\Delta v_{j}+1 \leq Ct_{j}^{-2}.
\end{equation}
Combining this with Theorem \ref{global estimate},
\[
\Delta v_{j}(z) \leq C\exp\left(\frac{C\log M_{R_{j}/2}}{R_{j}^{2}+1}\cdot v_{j}(z)\right) \ \ \text{in $B_{R_{j}/2}$}.
\]
Fixing $L\leq R_{j}/6$ and restricting the above on $B_{3L}$, together with Lemma \ref{uniform quadratic growth}, we obtain
\[
\Delta v_{j}(z) \leq C\exp\left(\frac{C\log M_{R_{j}/2}}{R_{j}^{2}+1}\cdot C(L^{2}+1)\right) \ \ \text{in $B_{3L}$}.
\]
Combining this with \eqref{M R j 2 estimate},
\begin{equation}\label{uniform Laplacian estimate}
\Delta v_{j}(z) \leq C\exp\left(\frac{C\log(Ct_{j}^{-2})}{R_{j}^{2}+1}\cdot C(L^{2}+1)\right) \ \ \text{in $B_{3L}$}.
\end{equation}
Since $\log(t_{j}^{-2})=o(R_{j}^{2})$, then the above implies $\Delta v_{j}\leq C_{L}$ in $B_{3L}$. Using \cite[Theorem 1]{RS84}, we have $\|\nabla\Delta v_{j}\|_{B_{2L}}\leq C_{L}$.
Applying the Schauder estimate and bootstrapping argument, we obtain $\|v_{j}\|_{C^{k}(B_{L})}\leq C_{k,L}$ for any positive integer $k$. Passing to a subsequence, there exists a function $Q\in\PSH(\C^{n})\cap C^{\infty}(\C^{n})$ such that $v_{j}\to Q$ in $C_{\mathrm{loc}}^{\infty}(\C^{n})$.

Next we will show that $Q$ is a quadratic polynomial. It follows from Lemma \ref{uniform quadratic growth} that
\[
C^{-1}(|z|^{2}+1) \leq Q(z) \leq C(|z|^{2}+1) \ \ \text{on $\mathbb{C}^{n}$}.
\]
Using \eqref{uniform Laplacian estimate} and $\log(t_{j}^{-2})=o(R_{j}^{2})$, we obtain $\Delta Q\leq C$ on $\C^{n}$. By \cite[Theorem 2]{RS84}, $Q_{i\ov{j}}$ are constant. Set
\[
H(z) := Q(z)-Q_{i\ov{j}}(0)z^{i}\ov{z}^{j}.
\]
Then $H$ is a harmonic function on $\C^{n}$ with quadratic growth. By the classical Liouville theorem, $H$ must be a quadratic polynomial, which implies $Q$ is also a quadratic polynomial, as required.
\end{proof}

\end{document}